\documentclass{amsart}
\usepackage{amssymb}
\usepackage{amsthm}
 
\theoremstyle{definition} 
\newtheorem{ex}{Example}[section] \theoremstyle{remark}
\newtheorem{rem}{Remark}[section] \newcommand{\pn}{\par\noindent}
\newcommand{\pmn}{\par\medskip\noindent}
\newcommand{\pbn}{\par\bigskip\noindent}
\begin{document}
\title{A representation of a set of maps as a ribbon bipartite graph}
\author{Yury Kochetkov}
\date{}
\begin{abstract} In this purely experimental work we try to represent the set
of plane maps with 3 vertices and 3 faces as a bipartite ribbon graph.
In particular, this construction allows one to estimate the genus of
the initial set.\end{abstract}
\email{yukochetkov@hse.ru}
\maketitle

\section{Introduction}
\pn By $M_{3,3}$ will be denoted the set of plane maps with 3
vertices and 3 faces (we chose this set because our reasoning is
nontrivial on one hand, but not very complex on another). \pmn A
set of cubic maps can be represented as a graph, where two maps
$A$ and $B$ are adjacent, if there exists a flip, that transforms
map $A$ into map $B$ (see \cite{BH}):
\[\begin{picture}(200,70) \multiput(0,15)(0,40){2}{\line(1,0){10}}
\multiput(10,5)(0,50){2}{\line(0,1){10}}
\multiput(10,15)(10,30){2}{\line(1,1){10}}
\multiput(10,55)(10,-30){2}{\line(1,-1){10}}
\multiput(30,15)(0,40){2}{\line(1,0){10}}
\multiput(30,5)(0,50){2}{\line(0,1){10}}
\put(20,25){\line(0,1){20}}

\put(65,32){$\Rightarrow$}

\multiput(100,25)(0,20){2}{\line(1,0){10}}
\multiput(110,15)(0,30){2}{\line(0,1){10}}
\multiput(110,25)(30,10){2}{\line(1,1){10}}
\multiput(110,45)(30,-10){2}{\line(1,-1){10}}
\multiput(150,25)(0,20){2}{\line(1,0){10}}
\multiput(150,15)(0,30){2}{\line(0,1){10}}
\put(120,35){\line(1,0){20}}
\end{picture}\]
\begin{center}{Figure 1. A flip}\end{center}
\pmn In our case there are no flips, so we had to define when two
maps from $M_{3,3}$ are "adjacent". \pmn We cannot transform a
given map into a triangulation, but can --- into a
\emph{quadrangulation}, i.e. into a bipartite map, where all faces
have perimeter 4 (see \cite{CM}, \cite{CV}). The transformation is
performed in four steps: a) paint white all vertices of the given
map; b) into each face ("a country") put a black vertex (the
"capital"); c) inside each county connect its capital with white
vertices in its border; d) delete initial edges.
\begin{ex}
\[\begin{picture}(400,80) \put(18,40){\oval(40,40)[l]}
\put(62,40){\oval(40,40)[r]}
\multiput(22,20)(0,20){3}{\line(1,0){36}}
\multiput(20,22)(0,20){2}{\line(0,1){16}}
\multiput(60,22)(0,20){2}{\line(0,1){16}}
\multiput(20,20)(0,20){3}{\circle{4}}
\multiput(60,20)(0,20){3}{\circle{4}}

\put(158,40){\oval(40,40)[l]} \put(202,40){\oval(40,40)[r]}
\multiput(162,20)(0,20){3}{\line(1,0){36}}
\multiput(160,22)(0,20){2}{\line(0,1){16}}
\multiput(200,22)(0,20){2}{\line(0,1){16}}
\multiput(160,20)(0,20){3}{\circle{4}}
\multiput(200,20)(0,20){3}{\circle{4}}
\multiput(145,40)(70,0){2}{\circle*{3}}
\multiput(180,30)(0,20){2}{\circle*{3}} \put(240,40){\circle*{3}}
\qbezier[20](145,40)(152,50)(159,59)
\qbezier[20](145,40)(152,30)(159,21)
\qbezier[15](145,40)(151,40)(158,40)
\qbezier[20](215,40)(208,50)(201,59)
\qbezier[20](215,40)(208,30)(201,21)
\qbezier[15](215,40)(208,40)(202,40)
\qbezier[40](161,41)(180,50)(199,59)
\qbezier[40](161,59)(180,50)(199,41)
\qbezier[40](161,21)(180,30)(199,39)
\qbezier[40](161,39)(180,30)(199,21)

\qbezier[60](240,40)(230,70)(201,61)
\qbezier[60](240,40)(230,10)(201,19)
\qbezier[80](240,40)(240,100)(161,61)
\qbezier[80](240,40)(240,-20)(161,19)

\put(105,37){$\Rightarrow$}

\put(308,40){\oval(36,40)[l]} \put(352,40){\oval(36,40)[r]}
\multiput(310,20)(0,20){3}{\circle{4}}
\multiput(350,20)(0,20){3}{\circle{4}}
\multiput(290,40)(80,0){2}{\circle*{3}}
\multiput(330,30)(0,20){2}{\circle*{3}}
\multiput(312,21)(0,20){2}{\line(2,1){37}}
\multiput(312,39)(0,20){2}{\line(2,-1){37}}
\multiput(290,40)(62,0){2}{\line(1,0){18}}
\put(390,40){\circle*{3}}

\qbezier(390,40)(380,70)(351,61) \qbezier(390,40)(380,10)(351,19)
\qbezier(390,40)(390,100)(311,61)
\qbezier(390,40)(390,-20)(311,19) \put(260,37){$\Rightarrow$}
\end{picture}\]
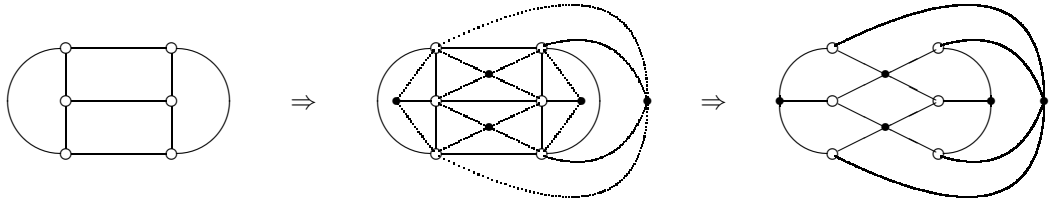
\begin{center}{Figure 2. From map to quadrangulation}\end{center}
\end{ex}
\pn
\begin{rem} It must be noted that the transformation "map $\Rightarrow$
quadrangulation" is invertible. \end{rem}
\pn
We want to construct a bipartite graph $G$, where white vertices are
$M_{3,3}$ maps and black vertices are $M(4,4,6)$ maps, i.e. bipartite
maps with 6 vertices and 3 faces: two of the perimeter 4 and one of the
perimeter 6. A map from $M(4,4,6)$ generates three quadrangulations.
Here are examples.
\begin{ex}
\[\begin{picture}(340,55) \multiput(0,20)(40,0){2}{\circle*{3}}
\multiput(0,40)(40,0){2}{\circle{4}} \put(20,20){\circle{4}}
\put(20,40){\circle*{3}} \multiput(0,20)(20,20){2}{\line(1,0){18}}
\multiput(0,20)(40,0){2}{\circle*{3}}
\multiput(2,40)(20,-20){2}{\line(1,0){18}}
\multiput(0,20)(40,0){2}{\circle*{3}}
\multiput(0,20)(40,0){2}{\line(0,1){18}}
\put(20,22){\line(0,1){18}}

\put(60,27){$\Rightarrow$}

\multiput(90,20)(40,0){2}{\circle*{3}}
\multiput(90,40)(40,0){2}{\circle{4}} \put(110,20){\circle{4}}
\put(110,40){\circle*{3}}
\multiput(90,20)(20,20){2}{\line(1,0){18}}
\multiput(90,20)(40,0){2}{\circle*{3}}
\multiput(92,40)(20,-20){2}{\line(1,0){18}}
\multiput(90,20)(40,0){2}{\circle*{3}}
\multiput(90,20)(40,0){2}{\line(0,1){18}}
\put(110,22){\line(0,1){18}} \put(130,40){\oval(40,20)[t]}
\put(130,18){\oval(40,18)[b]} \put(150,18){\line(0,1){22}}

\put(170,27){\small or}

\multiput(195,20)(40,0){2}{\circle*{3}}
\multiput(195,40)(40,0){2}{\circle{4}} \put(215,20){\circle{4}}
\put(215,40){\circle*{3}}
\multiput(195,20)(20,20){2}{\line(1,0){18}}
\multiput(195,20)(40,0){2}{\circle*{3}}
\multiput(197,40)(20,-20){2}{\line(1,0){18}}
\multiput(195,20)(40,0){2}{\circle*{3}}
\multiput(195,20)(40,0){2}{\line(0,1){18}}
\put(215,22){\line(0,1){18}} \put(220,20){\oval(50,22)[b]}
\put(237,20){\oval(16,42)[tr]}

\put(265,27){\small or}

\multiput(290,20)(40,0){2}{\circle*{3}}
\multiput(290,40)(40,0){2}{\circle{4}} \put(310,20){\circle{4}}
\put(310,40){\circle*{3}}
\multiput(290,20)(20,20){2}{\line(1,0){18}}
\multiput(290,20)(40,0){2}{\circle*{3}}
\multiput(292,40)(20,-20){2}{\line(1,0){18}}
\multiput(290,20)(40,0){2}{\circle*{3}}
\multiput(290,20)(40,0){2}{\line(0,1){18}}
\put(310,22){\line(0,1){18}} \put(315,42){\oval(50,22)[t]}
\put(330,42){\oval(20,44)[br]}
\end{picture}\]
\[\begin{picture}(340,50) \multiput(0,25)(15,0){2}{\circle*{3}}
\put(30,25){\circle{4}} \qbezier(0,26)(15,45)(29,27)
\qbezier(0,24)(15,5)(29,23) \put(15,26){\line(1,0){13}}
\put(28,25){\oval(56,40)[l]} \multiput(30,5)(0,40){2}{\circle{4}}
\put(32,25){\oval(46,40)[r]} \put(55,25){\circle*{3}}

\put(70,23){$\Rightarrow$}

\multiput(95,25)(15,0){2}{\circle*{3}} \put(125,25){\circle{4}}
\qbezier(95,26)(110,45)(124,27) \qbezier(95,24)(110,5)(124,23)
\put(110,26){\line(1,0){13}} \put(123,25){\oval(56,40)[l]}
\multiput(125,5)(0,40){2}{\circle{4}}
\put(127,25){\oval(46,40)[r]} \put(150,25){\circle*{3}}
\put(127,26){\line(1,0){23}}

\put(165,23){\small or}

\multiput(190,25)(15,0){2}{\circle*{3}} \put(220,25){\circle{4}}
\qbezier(190,26)(205,45)(219,27) \qbezier(190,24)(205,5)(219,23)
\put(205,26){\line(1,0){13}} \put(218,25){\oval(56,40)[l]}
\multiput(220,5)(0,40){2}{\circle{4}}
\put(222,25){\oval(46,40)[r]} \put(245,25){\circle*{3}}
\qbezier(190,25)(195,40)(210,40) \qbezier(210,40)(230,40)(230,25)
\qbezier(230,25)(230,15)(221,7)

\put(260,23){\small or}

\multiput(285,25)(15,0){2}{\circle*{3}} \put(315,25){\circle{4}}
\qbezier(285,26)(300,45)(314,27) \qbezier(285,24)(300,5)(314,23)
\put(300,26){\line(1,0){13}} \put(313,25){\oval(56,40)[l]}
\multiput(315,5)(0,40){2}{\circle{4}}
\put(317,25){\oval(46,40)[r]} \put(340,25){\circle*{3}}
\qbezier(285,25)(290,10)(305,10) \qbezier(305,10)(325,10)(325,25)
\qbezier(325,25)(325,35)(316,43)
\end{picture}\]
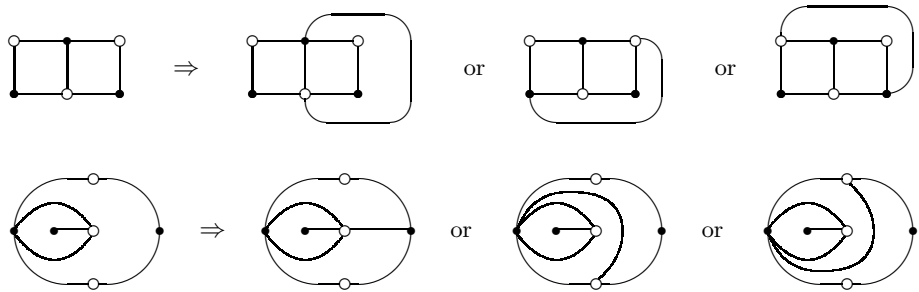
\begin{center}{Figure 3. Generation of quadrangulations}\end{center}
\end{ex}
\pn I.e. we connect a black vertex in the border of the face of
perimeter 6 with the opposite white vertex by an arc inside this
face. Clearly there are three possibilities to do it. Now the
construction of a bipartite graph is clear: a black vertex (i.e.
some $M(4,4,6)$ map $M$) is adjacent to a white vertex (i.e. some
$M_{3,3}$ map $N$) only when $M$ generates a quadrangulation of
$N$. Thus, the $M(4,4,6)$ map
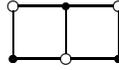
\[\begin{picture}(40,30) \multiput(0,5)(40,0){2}{\circle*{3}}
\multiput(0,25)(40,0){2}{\circle{4}} \put(20,5){\circle{4}}
\put(20,25){\circle*{3}} \multiput(0,5)(20,20){2}{\line(1,0){18}}
\multiput(0,5)(40,0){2}{\circle*{3}}
\multiput(2,25)(20,-20){2}{\line(1,0){18}}
\multiput(0,5)(40,0){2}{\circle*{3}}
\multiput(0,5)(40,0){2}{\line(0,1){18}} \put(20,7){\line(0,1){18}}
\end{picture}\] is connected to the map
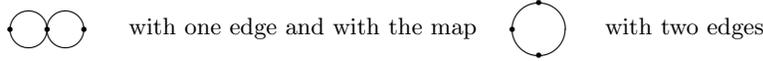
\[\begin{picture}(285,30) \multiput(0,15)(14,0){3}{\circle*{2}}
\multiput(7,15)(14,0){2}{\circle{14}}

\put(45,13){\small with one edge and with the map}

\put(200,15){\circle{20}} \put(190,15){\circle*{2}}
\multiput(200,5)(0,20){2}{\circle*{2}}

\put(225,13){\small with two edges}
\end{picture}\] However, we want to construct not a bipartite
graph, but bipartite \emph{ribbon} graph, i.e. we have to define a
cyclic order of outgoing edges in black and white vertices. It is
easy to do for black vertices. Given a map $M\in M(4,4,6)$, we
make a circuit of its face of perimeter 6, having the face to the
right of our path. This circuit gives us a cyclic order of black
vertices.
\begin{ex} For the $M(4,4,6)$ map
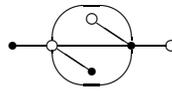
\[\begin{picture}(60,40) \multiput(0,20)(45,0){2}{\circle*{3}}
\multiput(15,20)(45,0){2}{\circle{4}} \put(0,20){\line(1,0){13}}
\put(17,20){\line(1,0){30}} \put(45,20){\line(1,0){13}}
\put(30,22){\oval(30,26)[t]} \put(30,18){\oval(30,26)[b]}
\put(30,10){\circle*{3}} \put(30,10){\line(-3,2){13}}
\put(30,30){\circle{4}} \put(45,20){\line(-3,2){13}}
\end{picture}\] we have the following cyclic order of $M_{3,3}$ maps,
adjacent to it:
\[\begin{picture}(180,30) \qbezier(0,15)(0,30)(20,15)
\qbezier(0,15)(0,0)(20,15) \qbezier(40,15)(40,30)(20,15)
\qbezier(40,15)(40,0)(20,15) \put(20,15){\circle*{3}}
\put(20,15){\line(-1,2){6}} \put(20,15){\line(1,2){6}}
\multiput(14,27)(12,0){2}{\circle*{2}} \put(19,3){\tiny 1}

\qbezier(70,15)(70,30)(90,15) \qbezier(70,15)(70,0)(90,15)
\qbezier(110,15)(110,30)(90,15) \qbezier(110,15)(110,0)(90,15)
\put(90,15){\circle*{3}} \put(90,5){\line(0,1){20}}
\multiput(90,5)(0,20){2}{\circle*{2}} \put(93,3){\tiny 2}

\qbezier(140,15)(140,30)(160,15) \qbezier(140,15)(140,0)(160,15)
\qbezier(180,15)(180,30)(160,15) \qbezier(180,15)(180,0)(160,15)
\put(160,15){\circle*{3}} \put(140,15){\circle*{2}}
\put(160,27){\circle*{2}} \put(160,15){\line(0,1){12}}
\put(159,3){\tiny 3} \end{picture}\]
\begin{center}{Figure 4. Cyclic order}\end{center}
\end{ex}
\pn It is more difficult to define an analogous cyclic order at
white vertices of $G$. It can be done in the following way. Let
$M$ be a $M_{3,3}$ map and $Q$ --- its quadrangulation. If we
delete an edge of $Q$ that separates two faces, then we obtain a
$M(4,4,6)$ map. Thus, it is enough to define a cyclic order in the
set of those edges of $Q$ that separate two faces. Let $Q^*$ be
the dual map. An edge in $Q$ that does not separate two faces
generates a loop in $Q^*$. As all vertices of $Q^*$ have degree 4,
then it is an Eulerian graph and has a cycle that contains all
edges except loops. Faces of the map $Q^*$ can be painted in black
and white in the chess like mode. Let external face be black. In
moving along the cycle we will have black faces to the right (this
agreement do not define a cycle uniquely).
\begin{ex}
\[\begin{picture}(350,40) \put(15,25){\circle{30}}
\put(15,10){\line(0,1){30}} \multiput(15,10)(0,30){2}{\circle*{2}}
\put(31,25){\circle*{2}} \put(13,1){\tiny M}

\multiput(60,25)(30,0){2}{\circle*{3}}
\multiput(90,10)(0,30){2}{\circle{4}} \put(75,25){\circle{4}}
\put(105,25){\circle*{3}} \put(60,25){\line(1,0){13}}
\put(60,25){\line(2,1){28}} \put(60,25){\line(2,-1){28}}
\put(77,25){\line(1,0){13}} \put(90,12){\line(0,1){26}}
\put(105,25){\line(-1,1){14}} \put(105,25){\line(-1,-1){14}}
\put(88,1){\tiny Q}

\put(155,25){\oval(40,20)} \put(155,25){\oval(40,30)}
\put(155,25){\oval(10,20)} \multiput(135,25)(40,0){2}{\circle*{2}}
\multiput(155,15)(0,20){2}{\circle*{2}} \put(153,1){\tiny $Q^*$}
\put(178,23){\tiny A} \put(128,23){\tiny B} \put(153,29){\tiny C}
\put(153,17){\tiny D}

\put(210,23){\small ABCDBADCA} \put(215,1){\footnotesize Eulerian
cycle}

\multiput(300,25)(30,0){2}{\circle*{3}}
\multiput(330,10)(0,30){2}{\circle{4}} \put(315,25){\circle{4}}
\put(345,25){\circle*{3}} \put(300,25){\line(1,0){13}}
\put(300,25){\line(2,1){28}} \put(300,25){\line(2,-1){28}}
\put(317,25){\line(1,0){13}} \put(330,12){\line(0,1){26}}
\put(345,25){\line(-1,1){14}} \put(345,25){\line(-1,-1){14}}

\put(337,34){\scriptsize 1} \put(313,34){\scriptsize 2}
\put(308,24){\tiny 3} \put(313,11){\scriptsize 4}
\put(337,11){\scriptsize 5} \put(326,17){\tiny 6}
\put(320,26){\tiny 7} \put(326,30){\tiny 8} \put(295,1){\tiny
Order of edges in Q}
\end{picture}\]
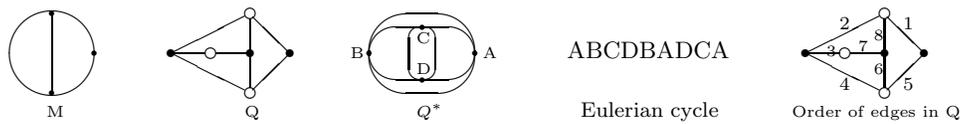
\begin{center}{Figure 5. Cyclic order of edges at a white
vertex}\end{center}
\end{ex}

\section{$M_{3,3}$ maps}
\pn Here we enumerate all $M_{3,3}$ maps. For the first two maps
we also present quadrangulations, dual to quadrangulations and the
order of edges. Loops are deleted.
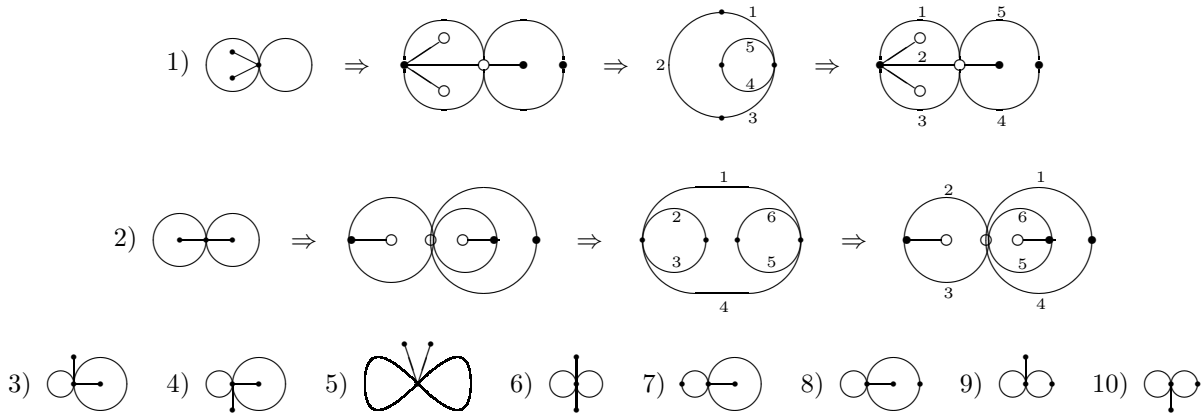
\[\begin{picture}(330,60) \put(0,28){1)}
\multiput(25,30)(20,0){2}{\circle{20}}
\multiput(25,25)(0,10){2}{\circle*{2}} \put(25,25){\line(2,1){10}}
\put(25,35){\line(2,-1){10}} \put(35,30){\circle*{2}}

\put(67,27){$\Rightarrow$}

\put(90,30){\circle*{3}} \put(120,30){\circle{4}}
\put(105,32){\oval(30,30)[t]} \put(105,28){\oval(30,30)[b]}
\put(90,30){\line(1,0){28}} \multiput(105,20)(0,20){2}{\circle{4}}
\put(90,30){\line(3,2){13}} \put(90,30){\line(3,-2){13}}
\multiput(135,30)(15,0){2}{\circle*{3}}
\put(135,32){\oval(30,30)[t]} \put(135,28){\oval(30,30)[b]}
\put(122,30){\line(1,0){13}}

\put(165,27){$\Rightarrow$}

\put(210,30){\circle{40}} \put(220,30){\circle{20}}
\multiput(210,10)(0,40){2}{\circle*{2}}
\multiput(210,30)(20,0){2}{\circle*{2}} \put(220,48){\tiny 1}
\put(185,28){\tiny 2} \put(220,8){\tiny 3} \put(219,21){\tiny 4}
\put(219,35){\tiny 5}

\put(245,27){$\Rightarrow$}

\put(270,30){\circle*{3}} \put(300,30){\circle{4}}
\put(285,32){\oval(30,30)[t]} \put(285,28){\oval(30,30)[b]}
\put(270,30){\line(1,0){28}}
\multiput(285,20)(0,20){2}{\circle{4}}
\put(270,30){\line(3,2){13}} \put(270,30){\line(3,-2){13}}
\multiput(315,30)(15,0){2}{\circle*{3}}
\put(315,32){\oval(30,30)[t]} \put(315,28){\oval(30,30)[b]}
\put(302,30){\line(1,0){13}} \put(284,48){\tiny 1}
\put(284,31){\tiny 2} \put(284,7){\tiny 3} \put(314,7){\tiny 4}
\put(314,48){\tiny 5}
\end{picture}\]
\[\begin{picture}(370,60) \put(0,28){2)}
\multiput(25,30)(20,0){2}{\circle{20}} \put(25,30){\line(1,0){20}}
\multiput(25,30)(10,0){3}{\circle*{2}}

\put(67,27){$\Rightarrow$}

\put(105,30){\circle{30}} \put(133,30){\circle{24}}
\put(90,30){\circle*{3}} \multiput(105,30)(15,0){2}{\circle{4}}
\put(90,30){\line(1,0){13}} \put(132,30){\circle{4}}
\put(144,30){\circle*{3}} \put(134,30){\line(1,0){12}}
\put(140,30){\circle{40}} \put(160,30){\circle*{3}}

\put(175,27){$\Rightarrow$}

\put(230,30){\oval(60,40)} \multiput(212,30)(36,0){2}{\circle{24}}
\multiput(200,30)(24,0){2}{\circle*{2}}
\multiput(236,30)(24,0){2}{\circle*{2}} \put(229,52){\tiny 1}
\put(229,3){\tiny 4} \put(211,37){\tiny 2} \put(211,20){\tiny 3}
\put(247,20){\tiny 5} \put(247,37){\tiny 6}

\put(275,27){$\Rightarrow$}

\put(315,30){\circle{30}} \put(343,30){\circle{24}}
\put(300,30){\circle*{3}} \multiput(315,30)(15,0){2}{\circle{4}}
\put(300,30){\line(1,0){13}} \put(342,30){\circle{4}}
\put(354,30){\circle*{3}} \put(344,30){\line(1,0){12}}
\put(350,30){\circle{40}} \put(370,30){\circle*{3}}
\put(349,52){\tiny 1} \put(314,47){\tiny 2} \put(314,8){\tiny 3}
\put(349,4){\tiny 4} \put(342,19){\tiny 5}  \put(342,37){\tiny 6}
\end{picture}\]
\[\begin{picture}(450,40) \put(0,17){3)} \put(20,20){\circle{10}}
\put(35,20){\circle{20}} \put(25,20){\line(1,0){10}}
\put(25,20){\line(0,1){10}} \multiput(25,20)(10,0){2}{\circle*{2}}
\put(25,30){\circle*{2}}

\put(60,17){4)} \put(80,20){\circle{10}} \put(95,20){\circle{20}}
\put(85,20){\line(1,0){10}} \put(85,20){\line(0,-1){10}}
\multiput(85,20)(10,0){2}{\circle*{2}} \put(85,10){\circle*{2}}

\put(120,17){5)} \qbezier(135,20)(135,40)(155,20)
\qbezier(135,20)(135,0)(155,20) \qbezier(155,20)(175,40)(175,20)
\qbezier(155,20)(175,0)(175,20) \put(155,20){\circle*{2}}
\put(155,20){\line(1,3){5}} \put(155,20){\line(-1,3){5}}
\multiput(150,35)(10,0){2}{\circle*{2}}

\put(190,17){6)} \multiput(210,20)(10,0){2}{\circle{10}}
\put(215,10){\line(0,1){20}}
\multiput(215,10)(0,10){3}{\circle*{2}}

\put(240,17){7)} \put(260,20){\circle{10}}
\put(275,20){\circle{20}} \put(265,20){\line(1,0){10}}
\multiput(255,20)(10,0){3}{\circle*{2}}

\put(300,17){8)} \put(320,20){\circle{10}}
\put(335,20){\circle{20}} \put(325,20){\line(1,0){10}}
\multiput(325,20)(10,0){3}{\circle*{2}}

\put(360,17){9)} \multiput(380,20)(10,0){2}{\circle{10}}
\put(385,20){\line(0,1){10}}
\multiput(385,20)(10,0){2}{\circle*{2}} \put(385,30){\circle*{2}}

\put(410,17){10)} \multiput(435,20)(10,0){2}{\circle{10}}
\put(440,20){\line(0,-1){10}} \put(440,10){\circle*{2}}
\multiput(440,20)(10,0){2}{\circle*{2}}
\end{picture}\]
\[\begin{picture}(440,40)
\put(0,17){11)} \put(25,20){\circle{10}} \put(45,20){\oval(30,20)}
\put(30,20){\line(1,0){20}} \multiput(30,20)(10,0){3}{\circle*{2}}

\put(75,17){12)} \multiput(100,20)(10,0){2}{\circle{10}}
\put(105,20){\line(0,1){16}}
\multiput(105,20)(0,8){3}{\circle*{2}}

\put(130,17){13)} \put(155,20){\circle{10}}
\put(150,20){\line(1,0){20}}
\multiput(150,20)(10,0){3}{\circle*{2}}

\put(185,17){14)} \multiput(210,20)(10,0){2}{\circle{10}}
\put(225,20){\line(1,0){10}}
\multiput(215,20)(10,0){3}{\circle*{2}}

\put(250,17){15)} \put(275,20){\circle{10}}
\put(290,20){\circle{20}} \put(290,20){\line(1,0){10}}
\multiput(280,20)(10,0){3}{\circle*{2}}

\put(315,17){16)} \put(340,20){\circle{10}}
\put(365,20){\circle{20}} \put(345,20){\line(1,0){20}}
\multiput(345,20)(10,0){3}{\circle*{2}}

\put(390,17){17)} \multiput(415,20)(20,0){2}{\circle{10}}
\put(420,20){\line(1,0){10}} \put(420,20){\line(0,1){10}}
\multiput(420,20)(10,0){2}{\circle*{2}} \put(420,30){\circle*{2}}
\end{picture}\]
\[\begin{picture}(365,40)
\put(0,17){18)} \multiput(25,20)(20,0){2}{\circle{10}}
\put(30,20){\line(1,0){10}} \put(30,20){\line(0,-1){10}}
\multiput(30,20)(10,0){2}{\circle*{2}} \put(30,10){\circle*{2}}

\put(65,17){19)} \multiput(90,20)(10,0){2}{\circle{10}}
\multiput(85,20)(10,0){3}{\circle*{2}}

\put(120,17){20)} \put(150,20){\circle{20}}
\put(165,20){\circle{10}} \multiput(150,10)(0,20){2}{\circle*{2}}
\put(160,20){\circle*{2}}

\put(185,17){21)} \put(215,20){\circle{20}}
\put(205,20){\line(1,0){20}}
\multiput(205,20)(10,0){3}{\circle*{2}}

\put(240,17){22)} \multiput(265,20)(20,0){2}{\circle{10}}
\put(270,20){\line(1,0){10}}
\multiput(260,20)(10,0){3}{\circle*{2}}

\put(305,17){23)} \multiput(330,20)(30,0){2}{\circle{10}}
\multiput(335,20)(10,0){3}{\circle*{2}}
\put(335,20){\line(1,0){20}}
\end{picture}\]
\begin{center}{Figure 6. $M_{3,3}$ maps}\end{center}

\section{$M(4,4,6)$ maps}
\pn Here we enumerate $M(4,4,6)$ maps. Maps with numbers $n$ and $n'$ are
the same up to mirror symmetry or the swap of colors.
\[\begin{picture}(390,55) \put(15,32){\oval(30,16)[t]}
\put(15,28){\oval(30,16)[b]} \put(0,30){\circle*{3}}
\multiput(15,30)(15,0){2}{\circle{4}} \put(20,32){\oval(40,30)[t]}
\put(20,28){\oval(40,30)[b]} \put(40,30){\circle{4}}
\put(0,30){\line(1,0){13}} \multiput(52,18)(0,24){2}{\circle*{3}}
\put(41,31){\line(1,1){12}} \put(41,29){\line(1,-1){12}}
\put(20,2){\small $1$}

\put(85,32){\oval(30,16)[t]} \put(85,28){\oval(30,16)[b]}
\put(70,30){\circle{4}} \multiput(85,30)(15,0){2}{\circle*{3}}
\put(90,32){\oval(40,30)[t]} \put(90,28){\oval(40,30)[b]}
\put(110,30){\circle*{3}} \put(72,30){\line(1,0){13}}
\multiput(122,18)(0,24){2}{\circle{4}}
\put(110,30){\line(1,1){11}} \put(110,30){\line(1,-1){11}}
\put(90,2){\small $1'$}

\multiput(150,32)(20,0){2}{\oval(20,16)[t]}
\multiput(150,28)(20,0){2}{\oval(20,16)[b]}
\put(140,30){\circle*{3}} \multiput(150,30)(10,0){2}{\circle{4}}
\put(140,30){\line(1,0){8}}
\multiput(170,30)(10,0){2}{\circle*{3}}
\put(162,30){\line(1,0){8}} \put(190,30){\circle{4}}
\put(180,30){\line(1,0){8}} \put(160,2){\small $2$}

\multiput(220,32)(20,0){2}{\oval(20,16)[t]}
\multiput(220,28)(20,0){2}{\oval(20,16)[b]}
\put(210,30){\circle{4}} \multiput(220,30)(10,0){2}{\circle*{3}}
\put(212,30){\line(1,0){8}} \multiput(240,30)(10,0){2}{\circle{4}}
\put(230,30){\line(1,0){8}} \put(260,30){\circle*{3}}
\put(252,30){\line(1,0){8}} \put(230,2){\small $2'$}

\put(290,32){\oval(30,20)[t]} \put(290,28){\oval(30,20)[b]}
\multiput(275,30)(20,0){2}{\circle*{3}}
\multiput(285,30)(20,0){2}{\circle{4}} \put(275,30){\line(1,0){8}}
\put(295,30){\line(-1,0){8}} \put(295,32){\oval(40,30)[t]}
\put(295,28){\oval(40,30)[b]} \put(315,30){\circle{4}}
\put(317,30){\line(1,0){8}} \put(325,30){\circle*{3}}
\put(295,2){\small $3$}

\put(355,32){\oval(30,20)[t]} \put(355,28){\oval(30,20)[b]}
\multiput(340,30)(20,0){2}{\circle{4}}
\multiput(350,30)(20,0){2}{\circle*{3}}
\put(342,30){\line(1,0){8}} \put(358,30){\line(-1,0){8}}
\put(360,32){\oval(40,30)[t]} \put(360,28){\oval(40,30)[b]}
\put(380,30){\circle*{3}} \put(380,30){\line(1,0){8}}
\put(390,30){\circle{4}} \put(360,2){\small $3'$}
\end{picture}\]
\[\begin{picture}(360,55) \put(27,32){\oval(30,30)[t]}
\put(27,28){\oval(30,30)[b]} \put(12,30){\circle*{3}}
\put(42,30){\circle{4}} \put(12,30){\line(1,0){28}}
\multiput(30,22)(0,16){2}{\circle*{3}} \put(30,22){\line(3,2){11}}
\put(30,38){\line(3,-2){11}} \multiput(0,22)(0,16){2}{\circle{4}}
\put(12,30){\line(-3,2){10}} \put(12,30){\line(-3,-2){10}}
\put(25,2){\small $4$}

\put(85,32){\oval(30,30)[t]} \put(85,28){\oval(30,30)[b]}
\put(70,30){\circle{4}} \put(100,30){\circle*{3}}
\put(73,30){\line(1,0){28}} \multiput(58,22)(0,16){2}{\circle*{3}}
\put(58,22){\line(3,2){10}} \put(58,38){\line(3,-2){10}}
\multiput(88,22)(0,16){2}{\circle{4}}
\put(100,30){\line(-3,2){10}} \put(100,30){\line(-3,-2){10}}
\put(82,2){\small $4'$}

\put(115,30){\circle*{3}} \put(115,30){\line(1,0){8}}
\put(125,30){\circle{4}} \put(155,30){\circle*{3}}
\put(140,32){\oval(30,30)[t]} \put(140,28){\oval(30,30)[b]}
\put(127,30){\line(1,0){28}} \put(165,30){\circle{4}}
\put(155,30){\line(1,0){8}} \put(137,38){\circle*{3}}
\put(137,38){\line(-3,-2){10}} \put(143,22){\circle{4}}
\put(155,30){\line(-3,-2){10}} \put(138,2){\small $5$}

\put(180,30){\circle*{3}} \put(180,30){\line(1,0){8}}
\put(190,30){\circle{4}} \put(220,30){\circle*{3}}
\put(205,32){\oval(30,30)[t]} \put(205,28){\oval(30,30)[b]}
\put(192,30){\line(1,0){28}} \put(230,30){\circle{4}}
\put(220,30){\line(1,0){8}} \put(202,22){\circle*{3}}
\put(202,22){\line(-3,2){10}} \put(208,38){\circle{4}}
\put(220,30){\line(-3,2){10}} \put(203,2){\small $5'$}

\put(245,30){\circle{4}} \multiput(275,30)(20,0){2}{\circle*{3}}
\put(260,32){\oval(30,30)[t]} \put(260,28){\oval(30,30)[b]}
\put(285,30){\circle{4}} \put(247,30){\line(1,0){36}}
\put(287,30){\line(1,0){8}} \put(257,38){\circle*{3}}
\put(257,38){\line(-3,-2){10}} \put(263,22){\circle{4}}
\put(275,30){\line(-3,-2){10}} \put(260,2){\small $6$}

\put(310,30){\circle*{3}} \multiput(340,30)(20,0){2}{\circle{4}}
\put(325,32){\oval(30,30)[t]} \put(325,28){\oval(30,30)[b]}
\put(350,30){\circle*{3}} \put(310,30){\line(1,0){28}}
\put(342,30){\line(1,0){16}} \put(322,38){\circle{4}}
\put(310,30){\line(3,2){10}} \put(328,22){\circle*{3}}
\put(328,22){\line(3,2){10}} \put(325,2){\small $6'$}
\end{picture}\]
\[\begin{picture}(400,55) \multiput(0,30)(20,0){3}{\circle{4}}
\multiput(10,30)(20,0){3}{\circle*{3}}
\put(25,32){\oval(30,20)[t]} \put(25,28){\oval(30,20)[b]}
\multiput(2,30)(20,0){2}{\line(1,0){16}}
\put(42,30){\line(1,0){8}} \put(24,2){\small $7$}

\multiput(65,30)(20,0){3}{\circle{4}}
\multiput(75,30)(20,0){3}{\circle*{3}}
\put(80,32){\oval(30,20)[t]} \put(80,28){\oval(30,20)[b]}
\multiput(67,30)(20,0){2}{\line(1,0){16}}
\put(107,30){\line(1,0){8}} \put(79,2){\small $8$}

\multiput(130,30)(20,0){3}{\circle*{3}}
\multiput(140,30)(20,0){3}{\circle{4}}
\put(145,32){\oval(30,20)[t]} \put(145,28){\oval(30,20)[b]}
\multiput(142,30)(20,0){2}{\line(1,0){16}}
\put(130,30){\line(1,0){8}} \put(144,2){\small $8'$}

\put(195,30){\circle*{3}} \put(205,30){\circle{4}}
\put(195,30){\line(1,0){8}} \put(235,30){\circle*{3}}
\put(220,32){\oval(30,30)[t]} \put(207,30){\line(1,0){28}}
\put(223,38){\circle{4}} \put(235,30){\line(-3,2){10}}
\put(205,15){\circle*{3}} \put(235,15){\circle{4}}
\put(205,15){\line(0,1){13}} \put(205,15){\line(1,0){28}}
\put(235,17){\line(0,1){13}} \put(219,2){\small $9$}

\put(250,30){\circle{4}} \put(260,30){\circle*{3}}
\put(252,30){\line(1,0){8}} \put(290,30){\circle{4}}
\put(275,32){\oval(30,30)[t]} \put(260,30){\line(1,0){28}}
\put(278,38){\circle*{3}} \put(278,38){\line(3,-2){10}}
\put(260,15){\circle{4}} \put(290,15){\circle*{3}}
\put(260,17){\line(0,1){13}} \put(262,15){\line(1,0){28}}
\put(290,15){\line(0,1){13}} \put(274,2){\small $9'$}

\put(305,30){\circle{4}} \put(335,30){\circle*{3}}
\put(307,30){\line(1,0){28}} \put(320,32){\oval(30,30)[t]}
\put(317,38){\circle*{3}} \put(317,38){\line(-3,-2){10}}
\put(305,15){\circle*{3}} \put(335,15){\circle{4}}
\put(305,15){\line(0,1){13}} \put(305,15){\line(1,0){28}}
\put(335,17){\line(0,1){13}} \put(345,30){\circle{4}}
\put(335,30){\line(1,0){8}} \put(319,2){\small $10$}

\put(360,30){\circle*{3}} \put(390,30){\circle{4}}
\put(360,30){\line(1,0){28}} \put(375,32){\oval(30,30)[t]}
\put(372,38){\circle{4}} \put(360,30){\line(3,2){10}}
\put(360,15){\circle{4}} \put(390,15){\circle*{3}}
\put(360,17){\line(0,1){13}} \put(362,15){\line(1,0){28}}
\put(390,15){\line(0,1){13}} \put(400,30){\circle*{3}}
\put(392,30){\line(1,0){8}} \put(374,2){\small $10'$}
\end{picture}\]
\[\begin{picture}(425,55) \put(0,30){\circle{4}}
\multiput(10,30)(30,0){2}{\circle*{3}}
\multiput(25,15)(0,30){2}{\circle{4}} \put(23,30){\oval(26,30)[l]}
\put(27,30){\oval(26,30)[r]} \put(2,30){\line(1,0){8}}
\put(40,43){\oval(30,26)[lb]} \put(35,35){\circle*{3}}
\put(36,34){\line(-1,1){10}} \put(21,2){\small $11$}

\put(55,30){\circle*{3}} \multiput(65,30)(30,0){2}{\circle{4}}
\multiput(80,15)(0,30){2}{\circle*{3}}
\put(80,32){\oval(30,26)[t]} \put(80,28){\oval(30,26)[b]}
\put(55,30){\line(1,0){8}} \put(93,45){\oval(26,30)[lb]}
\put(90,35){\circle{4}} \put(79,46){\line(1,-1){10}}
\put(76,2){\small $11'$}

\multiput(110,30)(40,0){2}{\circle{4}}
\multiput(120,30)(10,0){3}{\circle*{3}}
\put(120,32){\oval(20,16)[t]} \put(120,28){\oval(20,16)[b]}
\put(140,32){\oval(20,16)[t]} \put(140,28){\oval(20,16)[b]}
\put(112,30){\line(1,0){8}} \put(140,30){\line(1,0){8}}
\put(130,45){\circle{4}} \put(130,30){\line(0,1){13}}
\put(128,2){\small $12$}

\multiput(165,30)(40,0){2}{\circle*{3}}
\multiput(175,30)(10,0){3}{\circle{4}}
\put(175,32){\oval(20,16)[t]} \put(175,28){\oval(20,16)[b]}
\put(195,32){\oval(20,16)[t]} \put(195,28){\oval(20,16)[b]}
\put(165,30){\line(1,0){8}} \put(197,30){\line(1,0){8}}
\put(185,45){\circle*{4}} \put(185,32){\line(0,1){13}}
\put(183,2){\small $12'$}

\put(220,30){\circle*{3}} \multiput(230,30)(10,0){2}{\circle{4}}
\put(230,32){\oval(20,20)[t]} \put(230,28){\oval(20,20)[b]}
\put(220,30){\line(1,0){8}}
\multiput(250,30)(20,0){2}{\circle*{3}} \put(260,30){\circle{4}}
\multiput(242,30)(20,0){2}{\line(1,0){8}}
\put(255,32){\oval(30,20)[t]} \put(255,28){\oval(30,20)[b]}
\put(245,2){\small $13$}

\put(285,30){\circle{4}} \multiput(295,30)(10,0){2}{\circle*{3}}
\put(295,32){\oval(20,20)[t]} \put(295,28){\oval(20,20)[b]}
\put(287,30){\line(1,0){8}} \multiput(315,30)(20,0){2}{\circle{4}}
\put(325,30){\circle*{3}}
\multiput(305,30)(20,0){2}{\line(1,0){8}}
\put(320,32){\oval(30,20)[t]} \put(320,28){\oval(30,20)[b]}
\put(310,2){\small $13'$}

\multiput(350,30)(20,0){2}{\circle{4}}
\multiput(370,20)(0,20){2}{\circle*{3}} \put(360,30){\circle*{3}}
\put(352,30){\line(1,0){8}} \put(370,20){\line(0,1){8}}
\put(370,40){\line(0,-1){8}} \put(370,32){\oval(40,16)[tl]}
\put(370,28){\oval(40,16)[bl]} \put(380,30){\circle{4}}
\put(370,32){\oval(20,16)[tr]} \put(370,28){\oval(20,16)[br]}
\put(365,2){\small $14$}

\multiput(395,30)(20,0){2}{\circle*{3}}
\multiput(415,20)(0,20){2}{\circle{4}} \put(405,30){\circle{4}}
\put(395,30){\line(1,0){8}} \put(415,22){\line(0,1){16}}
\put(413,30){\oval(36,20)[tl]} \put(413,30){\oval(36,20)[bl]}
\put(425,30){\circle*{3}} \put(417,30){\oval(16,20)[tr]}
\put(417,30){\oval(16,20)[br]} \put(410,2){\small $14'$}
\end{picture}\]
\[\begin{picture}(390,55) \multiput(0,30)(10,0){2}{\circle*{3}}
\put(20,30){\circle{4}} \put(10,32){\oval(20,16)[t]}
\put(10,28){\oval(20,16)[b]} \put(10,30){\line(1,0){8}}
\put(40,30){\circle*{3}} \multiput(20,10)(0,40){2}{\circle{4}}
\put(18,30){\oval(36,40)[l]} \put(22,30){\oval(36,40)[r]}
\put(17,0){\scriptsize $15$}

\multiput(55,30)(10,0){2}{\circle{4}} \put(75,30){\circle*{3}}
\put(65,32){\oval(20,16)[t]} \put(65,28){\oval(20,16)[b]}
\put(67,30){\line(1,0){8}} \put(95,30){\circle{4}}
\multiput(75,10)(0,40){2}{\circle*{3}}
\put(75,32){\oval(40,36)[t]} \put(75,28){\oval(40,36)[b]}
\put(72,0){\scriptsize $15'$}

\put(110,30){\circle*{3}} \multiput(110,20)(0,20){2}{\circle{4}}
\put(110,22){\line(0,1){16}} \put(130,30){\circle{4}}
\put(110,30){\line(1,0){18}}
\multiput(130,20)(0,20){2}{\circle*{3}}
\multiput(112,20)(0,20){2}{\line(1,0){18}}
\put(130,20){\line(0,1){8}} \put(130,32){\line(0,1){8}}
\put(117,2){\small $16$}

\multiput(145,30)(10,0){2}{\circle*{3}}
\multiput(165,30)(10,0){2}{\circle{4}}
\multiput(155,32)(20,0){2}{\oval(20,16)[t]}
\multiput(155,28)(20,0){2}{\oval(20,16)[b]}
\put(185,30){\circle*{3}} \put(195,30){\circle{4}}
\multiput(155,30)(30,0){2}{\line(1,0){8}}
\put(177,30){\line(1,0){8}} \put(163,2){\small $17$}

\multiput(210,30)(10,0){2}{\circle{4}}
\multiput(230,30)(10,0){2}{\circle*{3}}
\put(222,30){\line(1,0){8}}
\multiput(220,32)(20,0){2}{\oval(20,16)[t]}
\multiput(220,28)(20,0){2}{\oval(20,16)[b]}
\put(250,30){\circle{4}} \put(260,30){\circle*{3}}
\put(252,30){\line(1,0){8}} \put(228,2){\small $17'$}
\put(240,30){\line(1,0){8}}

\put(275,30){\circle*{3}} \multiput(305,30)(20,0){2}{\circle{4}}
\put(290,32){\oval(30,30)[t]} \put(290,28){\oval(30,30)[b]}
\put(275,30){\line(1,0){28}} \put(307,30){\line(1,0){16}}
\put(287,22){\circle{4}} \put(275,30){\line(3,-2){10}}
\put(293,38){\circle*{3}} \put(293,38){\line(3,-2){10}}
\put(315,30){\circle*{3}} \put(288,2){\small $18$}

\put(340,30){\circle{4}} \multiput(370,30)(20,0){2}{\circle*{3}}
\put(355,32){\oval(30,30)[t]} \put(355,28){\oval(30,30)[b]}
\put(342,30){\line(1,0){28}} \put(370,30){\line(1,0){8}}
\put(390,30){\line(-1,0){8}} \put(352,22){\circle*{3}}
\put(352,22){\line(-3,2){10}} \put(358,38){\circle{4}}
\put(370,30){\line(-3,2){10}} \put(380,30){\circle{4}}
\put(353,2){\small $18'$}
\end{picture}\]

\[\begin{picture}(380,55) \multiput(0,30)(20,0){2}{\circle*{3}}
\multiput(10,30)(20,0){2}{\circle{4}} \put(15,32){\oval(30,20)[t]}
\put(15,28){\oval(30,20)[b]} \put(40,30){\circle{4}}
\put(50,30){\circle*{3}} \put(40,32){\oval(20,16)[t]}
\put(40,28){\oval(20,16)[b]} \put(0,30){\line(1,0){8}}
\put(12,30){\line(1,0){8}} \put(42,30){\line(1,0){8}}
\put(22,2){\small $19$}

\multiput(65,30)(20,0){2}{\circle{4}}
\multiput(75,30)(20,0){2}{\circle*{3}}
\put(80,32){\oval(30,20)[t]} \put(80,28){\oval(30,20)[b]}
\put(105,30){\circle*{3}} \put(115,30){\circle{4}}
\put(105,32){\oval(20,16)[t]} \put(105,28){\oval(20,16)[b]}
\put(67,30){\line(1,0){16}} \put(105,30){\line(1,0){8}}
\put(87,2){\small $19'$}

\put(130,30){\circle{4}} \multiput(140,30)(10,0){2}{\circle*{3}}
\multiput(160,30)(10,0){2}{\circle{4}} \put(180,30){\circle*{3}}
\multiput(140,32)(30,0){2}{\oval(20,16)[t]}
\multiput(140,28)(30,0){2}{\oval(20,16)[b]}
\multiput(132,30)(40,0){2}{\line(1,0){8}}
\put(150,30){\line(1,0){8}} \put(151,2){\small $20$}

\put(215,30){\circle{4}} \multiput(195,30)(10,0){2}{\circle*{3}}
\multiput(235,30)(10,0){2}{\circle{4}} \put(225,30){\circle*{3}}
\multiput(205,32)(30,0){2}{\oval(20,16)[t]}
\multiput(205,28)(30,0){2}{\oval(20,16)[b]}
\multiput(205,30)(20,0){2}{\line(1,0){8}}
\put(217,30){\line(1,0){8}} \put(216,2){\small $21$}

\multiput(260,30)(40,0){2}{\circle{4}}
\multiput(280,13)(0,34){2}{\circle*{3}}
\put(280,32){\oval(40,30)[t]} \put(280,28){\oval(40,30)[b]}
\put(310,30){\circle*{3}} \put(270,35){\circle{4}}
\put(281,47){\line(-1,-1){10}} \put(302,30){\line(1,0){8}}
\qbezier(262,30)(285,25)(280,47) \put(277,2){\small $22$}

\multiput(320,30)(44,0){2}{\circle*{3}}
\multiput(342,13)(0,34){2}{\circle{4}}
\put(340,30){\oval(40,32)[l]} \put(344,30){\oval(40,32)[r]}
\put(374,30){\circle{4}} \put(330,35){\circle*{3}}
\put(330,35){\line(1,1){10}} \put(364,30){\line(1,0){8}}
\qbezier(320,30)(350,25)(342,45) \put(340,2){\small $22'$}
\end{picture}\]
\begin{center}{Figure 7. $M(4,4,6)$ maps}\end{center}

\section{The cyclic order of edges, outgoing from white vertices}
\pn In the table below in column denoted "N" numbers of $M_{3,3}$
maps are enumerated. In columns denoted "order" are enumerated (in
the cyclic order) numbers of $M(4,4,6)$ maps, generated by the
corresponding $M_{3,3}$ map.
\small{\[\begin{tabular}{|l|l|l|l|l|l|l|l|l|} \hline
\multicolumn{1}{|c|}{N} & \multicolumn{1}{c|}{order} &
\multicolumn{1}{c|}{N} & \multicolumn{1}{c|}{order} &
\multicolumn{1}{c|}{N} & \multicolumn{1}{c|}{order} &
\multicolumn{1}{c|}{N} & \multicolumn{1}{c|}{order} \\
\hline 1& $17,1',17,4',4'$& 7 & $15',8',8',15',19,13,3'$ & 13 &
$14',14',22,11',10',9,15'$ & 19 & $16,8,7,8',16,8,7,8'$\\ \hline

2 & $12',3',3',12',3',3'$ & 8 & $22',11,2,17,10',9$ & 14 &
$11',22,9',10,17',2'$ & 20 & $14,14,22',11,10,9',15$\\ \hline

3 & $12',13,2,6',6'$ & 9 & $6',5,6,9',7,9$ & 15 & $22,11',21,11,22',20$
&  21 & $16,14',14',16,16,16,14,14$\\ \hline

4 & $2,13,12',18,18$ & 10 & $10',7,10,18',5',18$ & 16 &
$20,19,19,20,19',19'$ & 22 & $8,8,15,19',13',3,15$ \\ \hline

5 & $5,4,5',4'$ & 11 & $21,1',1',21,1,1$ & 17 & $6,6,12,13',2'$ & 23 &
$12,3,3,12,3,3$\\ \hline

6 & $5',5,5',5$ & 12 & $17',1,17',4,4$ & 18 & $2',13',12,18',18'$ &&
 \\ \hline
\end{tabular}\]}
\begin{center}{Table 1}\end{center}

\section{The cyclic order of edges, outgoing from black vertices}
\pn
\[\begin{tabular}{|l|l|l|l|l|l|l|l|l|l|} \hline
\multicolumn{1}{|c|}{N} & \multicolumn{1}{c|}{order} &
\multicolumn{1}{c|}{N} & \multicolumn{1}{c|}{order} &
\multicolumn{1}{c|}{N} & \multicolumn{1}{c|}{order} &
\multicolumn{1}{c|}{N} & \multicolumn{1}{c|}{order} &
\multicolumn{1}{c|}{N} & \multicolumn{1}{c|}{order} \\ \hline
$1$&11,11,12 & $5$ & 9,6,5 & $9'$ & 20,14,9 & $13'$ & 22,17,18 &
$18$ & 10,4,4 \\ \hline $1'$ & 11,11,1 & $5'$ & 10,5,6 & $10$ &
20,10,14 & $14$ & 21,20,20 & $18'$ & 10,18,18 \\ \hline $2$ &
8,4,3 & $6$ & 17,17,9 & $10'$ & 13,10,8 & $14'$ & 21,13,13 & $19$ & 7,16,16 \\
\hline $2'$ & 14,18,17 & $6'$ & 3,3,9 & $11$ & 20,15,8 & $15$ &
20,22,22 & $19'$ & 22,16,16 \\ \hline $3$ & 22,23,23 & $7$ &
19,9,10 & $11'$ & 13,15,14 & $15'$ & 13,7,7 & $20$ & 15,16,16 \\
\hline $3'$ & 2,2,7 & $8$ & 19,22,22 & $12$ & 23,18,17 & $16$ &
19,21,21 & $21$ & 15,11,11 \\ \hline $4$ & 12,12,5 & $8'$ & 19,7,7
& $12'$ & 2,4,3 & $17$ & 8,1,1 & $22$ & 13,14,15 \\ \hline $4'$ &
1,1,5 & $9$ & 13,8,9 & $13$ & 7,3,4 & $17'$ & 14,12,12 & $22'$ & 20,8,15 \\
\hline
\end{tabular}\]
\begin{center}{Table 2}\end{center}

\section{The construction of a ribbon graph}
\pn We see, that our graph contains multiple edges. Thus, we have a multiple
edge that connects the first white vertex with the black vertex $4'$. Also
there are two multiple edges that connects the second white vertex
with the black vertex $3'$. However, edges that connect the first white
vertex with the black vertex $17$ \emph{do not} constitute a multiple edge.
In tables below we present cyclic orders in a simplified form, i.e. we
consider a multiple edge as an ordinary edge.
\[\begin{tabular}{|l|l|l|l|l|l|l|l|l|} \hline
\multicolumn{1}{|c|}{N} & \multicolumn{1}{c|}{order} &
\multicolumn{1}{c|}{N} & \multicolumn{1}{c|}{order} &
\multicolumn{1}{c|}{N} & \multicolumn{1}{c|}{order} &
\multicolumn{1}{c|}{N} & \multicolumn{1}{c|}{order} \\
\hline 1& $17,1',17,4'$& 7 & $15',8',15',19,13,3'$ & 13 &
$14',22,11',10',9,15'$ & 19 & $16,8,7,8'$\\ \hline

2 & $3',12'$ & 8 & $22',11,2,17,10',9$ & 14 &
$11',22,9',10,17',2'$ & 20 & $14,22',11,10,9',15$\\ \hline

3 & $12',13,2,6'$ & 9 & $6',5,6,9',7,9$ & 15 &
$22,11',21,11,22',20$ & 21 & $16,14',16,14$\\ \hline

4 & $2,13,12',18$ & 10 & $10',7,10,18',5',18$ & 16 &
$20,19,20,19'$ & 22 & $8,15,19',13',3,15$ \\ \hline

5 & $5,4,5',4'$ & 11 & $21,1',21,1$ & 17 & $6,12,13',2'$ & 23 & $3,12$\\
\hline

6 & $5',5$ & 12 & $17',1,17',4$ & 18 & $2',13',12,18'$ &&
\\ \hline \end{tabular}\]
\begin{center}{Table 3}\end{center}

\[\begin{tabular}{|l|l|l|l|l|l|l|l|l|l|} \hline
\multicolumn{1}{|c|}{N} & \multicolumn{1}{c|}{order} &
\multicolumn{1}{c|}{N} & \multicolumn{1}{c|}{order} &
\multicolumn{1}{c|}{N} & \multicolumn{1}{c|}{order} &
\multicolumn{1}{c|}{N} & \multicolumn{1}{c|}{order} &
\multicolumn{1}{c|}{N} & \multicolumn{1}{c|}{order} \\ \hline
$1$&11,12 & $5$ & 9,6,5 & $9'$ & 20,14,9 & $13'$ & 22,17,18 & $18$ & 10,4 \\
\hline $1'$ & 11,1 & $5'$ & 10,5,6 & $10$ & 20,10,14 & $14$ &
21,20 & $18'$ & 10,18 \\ \hline $2$ & 8,4,3 & $6$ & 17,9 & $10'$ &
13,10,8 & $14'$ & 21,13 & $19$ & 7,16 \\ \hline $2'$ & 14,18,17
& $6'$ & 3,9 & $11$ & 20,15,8 & $15$ & 20,22,22 & $19'$ & 22,16 \\
\hline $3$ & 22,23 & $7$ & 19,9,10 & $11'$ & 13,15,14 & $15'$ &
13,7,7 & $20$ & 15,16,16 \\ \hline $3'$ & 2,7 & $8$ & 19,22 & $12$
& 23,18,17 & $16$ & 19,21,21 & $21$ & 15,11,11 \\ \hline $4$ &
12,5 & $8'$ & 19,7 & $12'$ & 2,4,3 & $17$ & 8,1,1 & $22$ & 13,14,15 \\
\hline $4'$ & 1,5 & $9$ & 13,8,9 & $13$ & 7,3,4 & $17'$ & 14,12,12
& $22'$ & 20,8,15 \\ \hline
\end{tabular}\]
\begin{center}{Table 4}\end{center}
\pmn However, a structure of bipartite ribbon graph is not defined
yet. For example, when we come to the white vertex 1 from the
black vertex 17, we can go either to the black vertex $1'$, or to
the black vertex $4'$. All cases of such ambiguity in white
vertices are enumerated in the figure below:
\[\begin{picture}(365,40) \put(5,20){\circle{10}} \put(2,18){\tiny
$17$} \put(13,20){\vector(1,0){15}} \put(35,20){\circle{10}}
\put(34,18){\tiny 1} \multiput(65,10)(0,20){2}{\circle{10}}
\put(43,17){\vector(3,-1){14}} \put(43,23){\vector(3,1){14}}
\put(63,8){\tiny $4'$} \put(63,28){\tiny $1'$} \put(33,2){\small
a}

\put(100,20){\circle{12}} \put(96,18){\tiny $15'$}
\put(110,20){\vector(1,0){15}} \put(132,20){\circle{10}}
\put(130,18){\tiny 7} \multiput(162,10)(0,20){2}{\circle{10}}
\put(140,17){\vector(3,-1){14}} \put(140,23){\vector(3,1){14}}
\put(159,8){\tiny $19$} \put(160,28){\tiny $8'$}
\put(130,2){\small b}

\put(200,20){\circle{10}} \put(197,18){\tiny $22$}
\put(208,20){\vector(1,0){15}} \put(230,20){\circle{10}}
\put(228,18){\tiny 11} \multiput(260,10)(0,20){2}{\circle{10}}
\put(238,17){\vector(3,-1){14}} \put(238,23){\vector(3,1){14}}
\put(258,8){\tiny $1'$} \put(259,29){\tiny $1$} \put(228,2){\small
c}

\put(300,20){\circle{12}} \put(296,18){\tiny $17'$}
\put(310,20){\vector(1,0){15}} \put(332,20){\circle{10}}
\put(329,18){\tiny 12} \multiput(362,10)(0,20){2}{\circle{10}}
\put(340,17){\vector(3,-1){14}} \put(340,23){\vector(3,1){14}}
\put(361,8){\tiny $4$} \put(361,29){\tiny $1$} \put(330,2){\small
d}
\end{picture}\]
\[\begin{picture}(270,40) \put(5,20){\circle{10}} \put(2,18){\tiny
$21$} \put(13,20){\vector(1,0){15}} \put(35,20){\circle{10}}
\put(32,18){\tiny 16} \put(66,9){\circle{12}}
\put(65,30){\circle{10}} \put(43,17){\vector(3,-1){14}}
\put(43,23){\vector(3,1){14}} \put(62,7){\tiny $19'$}
\put(63,28){\tiny $19$} \put(33,2){\small e}

\put(100,20){\circle{10}} \put(98,18){\tiny $16$}
\put(108,20){\vector(1,0){15}} \put(130,20){\circle{10}}
\put(127,18){\tiny 21} \put(161,9){\circle{12}}
\put(160,30){\circle{10}} \put(138,17){\vector(3,-1){14}}
\put(138,23){\vector(3,1){14}} \put(157,7){\tiny $14'$}
\put(157,28){\tiny $14$} \put(128,2){\small f}

\put(200,20){\circle{10}} \put(197,18){\tiny $16$}
\put(208,20){\vector(1,0){15}} \put(230,20){\circle{10}}
\put(227,18){\tiny 22} \put(261,9){\circle{12}}
\put(260,30){\circle{10}} \put(238,17){\vector(3,-1){14}}
\put(238,23){\vector(3,1){14}} \put(257,7){\tiny $19'$}
\put(258,29){\tiny $8$} \put(228,2){\small n}
\end{picture}\]
\begin{center}{Figure 8}\end{center}
\pmn Analogously, in the figure below are enumerated all cases of
ambiguity in black vertices:
\[\begin{picture}(365,40) \put(5,20){\circle{10}} \put(2,18){\tiny
22} \put(13,20){\vector(1,0){15}} \put(35,20){\circle{10}}
\put(32,18){\tiny $15$} \multiput(65,10)(0,20){2}{\circle{10}}
\put(43,17){\vector(3,-1){14}} \put(43,23){\vector(3,1){14}}
\put(62,8){\tiny 22} \put(63,28){\tiny 20} \put(33,2){\small g}

\put(100,20){\circle{10}} \put(99,18){\tiny 7}
\put(108,20){\vector(1,0){15}} \put(131,20){\circle{12}}
\put(126,18){\tiny $15'$} \multiput(162,10)(0,20){2}{\circle{10}}
\put(140,17){\vector(3,-1){14}} \put(140,23){\vector(3,1){14}}
\put(159,8){\tiny 13} \put(161,28){\tiny 7} \put(130,2){\small h}

\put(200,20){\circle{10}} \put(197,18){\tiny 21}
\put(208,20){\vector(1,0){15}} \put(230,20){\circle{10}}
\put(227,18){\tiny $16$} \multiput(260,10)(0,20){2}{\circle{10}}
\put(238,16){\vector(3,-1){14}} \put(238,23){\vector(3,1){14}}
\put(257,8){\tiny 21} \put(257,28){\tiny 19} \put(229,2){\small i}

\put(300,20){\circle{10}} \put(299,18){\tiny 1}
\put(308,20){\vector(1,0){15}} \put(330,20){\circle{10}}
\put(327,18){\tiny $17$} \multiput(360,10)(0,20){2}{\circle{10}}
\put(338,17){\vector(3,-1){14}} \put(338,23){\vector(3,1){14}}
\put(359,8){\tiny 8} \put(359,29){\tiny 1} \put(330,2){\small j}
\end{picture}\]
\[\begin{picture}(270,40) \put(5,20){\circle{10}} \put(2,18){\tiny
12} \put(13,20){\vector(1,0){15}} \put(36,20){\circle{12}}
\put(31,18){\tiny $17'$} \multiput(65,10)(0,20){2}{\circle{10}}
\put(43,17){\vector(3,-1){14}} \put(43,23){\vector(3,1){14}}
\put(62,8){\tiny 14} \put(62,28){\tiny 12} \put(33,2){\small k}

\put(100,20){\circle{10}} \put(98,18){\tiny 16}
\put(108,20){\vector(1,0){15}} \put(130,20){\circle{10}}
\put(127,18){\tiny $21$} \multiput(160,10)(0,20){2}{\circle{10}}
\put(138,17){\vector(3,-1){14}} \put(138,23){\vector(3,1){14}}
\put(157,8){\tiny 16} \put(157,28){\tiny 15} \put(128,2){\small l}

\put(200,20){\circle{10}} \put(197,18){\tiny 11}
\put(208,20){\vector(1,0){15}} \put(230,20){\circle{10}}
\put(227,18){\tiny $22$} \multiput(260,10)(0,20){2}{\circle{10}}
\put(238,17){\vector(3,-1){14}} \put(238,23){\vector(3,1){14}}
\put(257,8){\tiny 15} \put(257,28){\tiny 11} \put(228,2){\small m}
\end{picture}\]
\begin{center}{Figure 9}\end{center}
\pmn If we come to, for example, the white vertex 1 from the black
vertex 17 for the first time, we have two options: upper outgoing
arrow (to the black vertex $1'$) --- option $a_1$, or lower
outgoing arrow (to the black vertex $4'$) --- option $a_2$.
Analogously, coming to the black vertex 21 from the white vertex
16 for the first time, we have two options: upper outgoing arrow
(to the white vertex 15) --- option $l_1$, or lower outgoing arrow
( to the white vertex 16) --- option $l_2$. \pmn If we want to
embed our graph into a compact surface we must choose an option in
each case ambiguity. So we have $2^{14}$ different variants of
embedding. However, a computer check demonstrates that there are
either 7 faces, or 9 faces, thus our graph is embedded either in a
surface of the genus 17, or in a surface of the genus 18.

\vspace{5mm}

\end{document}